\documentclass[11pt]{article}
\usepackage{amsfonts}
\usepackage{amscd}
\usepackage{amsmath}
\usepackage{epsfig}

\def\be{\begin{equation}}
\def\ee{\end{equation}}
\def\ben{\begin{displaymath}}
\def\een{\end{displaymath}}
\def\baa{\begin{eqnarray}}
\def\eaa{\end{eqnarray}}

\def\ba{\begin{array}}
\def\ea{\end{array}}

\makeatletter
\@addtoreset{equation}{section}
\makeatother

\usepackage{amsthm}
\usepackage{amsmath}
\usepackage{amssymb}
\usepackage{amsbsy}
\usepackage{amsfonts}
\usepackage{latexsym}
\usepackage{amscd}
\usepackage{eucal}
\usepackage{mathrsfs}
\usepackage[all, knot]{xy}
\xyoption{arc} \xyoption{color}\xyoption{line}

\newtheorem{prop}{Proposition}
\newtheorem{lem}{Lemma}

\newtheorem{Remark}{Remark}
\newtheorem{thm}{Theorem}
\newtheorem{defi}{Definition}

\begin{document}
\title{Flat conical Laplacian in the square of the canonical bundle and its regularized determinants}

\author{Alexey Kokotov}
\maketitle

\vskip2cm
{\bf Abstract.}
Let $X$ be a compact Riemann surface of genus $g\geq 2$ equipped with flat conical metric $|\Omega|$, where  $\Omega$ be a holomorphic quadratic differential on $X$ with $4g-4$ simple zeroes. Let $K$ be the canonical line bundle on $X$. Introduce the Cauchy-Riemann operators $\bar \partial$ and $\partial$  acting on sections of holomorphic line bundles over $X$  ($K^2$ in the definition of $\Delta^{(2)}_{|\Omega|}$ below)  and, respectively, anti-holomorphic line bundles ($\bar { K}^{-1}$ below). Consider the Laplace operator
$\Delta^{(2)}_{|\Omega|}:=|\Omega| \partial |\Omega|^{-2}\bar\partial$ acting in the Hilbert space of square integrable sections of the bundle $K^2$ equipped with inner product
$<Q_1, Q_2>_{K^2}=\int_X\frac {Q_1\bar Q_2}{|\Omega|}$.

We discuss two natural definitions of the determinant of the operator $\Delta^{(2)}_{|\Omega|}$. The first one uses the zeta-function of some special self-adjoint extension of the operator (initially defined on smooth sections of $K^2$ vanishing near the zeroes of $\Omega$), the second one is an analog of Eskin-Kontsevich-Zorich (EKZ) regularization of the determinant of the conical Laplacian acting in the trivial bundle.  In contrast to the situation of operators acting in the trivial bundle, for operators acting in $K^2$ these two regularizations turn out to be essentially different.  Considering the regularized determinant of $\Delta^{(2)}_{|\Omega|}$ as a functional on the moduli space $Q_g(1, \dots, 1)$ of quadratic differentials with simple zeroes on compact Riemann surfaces of genus $g$, we derive explicit expressions for this functional for the both regularizations. The expression for the EKZ regularization is closely related to the well-known explicit expressions for the Mumford measure on the moduli space of compact Riemann surfaces of genus $g$.

\vskip2cm

\section{Introduction}

Let $X$ be a compact Riemann surface of genus $g\geq 2$ and let $\Omega$ be a holomorphic quadratic differential with simple zeroes $P_1, \dots, P_{4g-4}$. Then its modulus, $|\Omega|$, gives a flat conformal metric on $X$ with conical singularities at the points $P_k$ with conical angles $3\pi$. The spectral theory of  scalar (i. e. acting in the trivial bundle over $X$)  Laplacians in flat conical metrics on Riemann surfaces is fairly well understood. In particular, an important spectral invariant, the $\zeta$-regularized determinant of the Laplacian, is explicitly computed in \cite{KokKor}  for metrics given by the modulus square of a holomorphic one-form and in \cite{KokKorLeipzig} for metrics given by the modulus of an arbitrary holomorphic quadratic differential; in \cite{KokProcAMS} these results are generalized for arbitrary flat conical metrics. In the present work we are dealing with Laplacian in metric $|\Omega|$ acting in the square $K^2$ of the canonical line bundle. Keeping in mind classical results on the explicit formulas for the determinant of the hyperbolic Laplacian in integer and half-integer powers of the canonical bundle (\cite{D'Hoker-Phong}, \cite{Sarnak}) we are going to define and then compute the regularized determinant of this operator. Our main motivation came from    paper \cite{Giddings}, where a similar object
was introduced and freely used by means of wonderful heuristics (see also \cite{Menotti}, where a closely related problem  for general flat conical metrics on Riemann surfaces of genera zero and one is considered).

The following immediate observation (essentially coinciding with the one from \cite{Giddings}, p. 106) is the starting point of the present work.
Let $$L_2(K^2; |\Omega|)=\{u\in \Gamma(K^2): \int_X|u|^2|\Omega|^{-1}< \infty\}$$ be the space of square integrable sections of the bundle $K^2$ and let $L_2(X; |\Omega|)=\{v\in \Gamma(E): \int_X|v|^2|\Omega|<\infty \}$ be the space of square integrable functions on $X$ (i. e. sections of the trivial bundle $E=X\times {\mathbb C}$). Then unbounded operators $$\Delta^{(2)}_{|\Omega|}:=|\Omega| \partial |\Omega|^{-2}\bar\partial: L_2(K^2; |\Omega|)\to
L_2(K^2; |\Omega|)$$ and $$\Delta^{(0)}_{|\Omega|}:= |\Omega|^{-1}\partial \bar \partial : L_2(X; |\Omega|)\to L_2(X; |\Omega|)$$ with domains consisting of smooth sections (functions) with supports on compact subsets of $X\setminus\{P_1, \dots, P_{4g-4}\}$ are unitary equivalent. The corresponding unitary equivalence is given by the operator $U: L_2(K^2; |\Omega|)\to L_2(X, |\Omega|)$ defined via
$$u \mapsto U(u)=u\Omega^{-1}\,.$$
 This could be the end of the story: the two operators are unitary equivalent, so their regularized determinants (as well as all other spectral invariants) are the same. Thus, everything is reduced to the known case of scalar conical Laplacians. Unfortunately, the situation is not that simple. It is reasonable to include into the domain of the operator $\Delta^{(2)}_\Omega$ its zero modes, holomorphic quadratic differentials on $X$. But if $q$ is a holomorphic quadratic differential on $X$ not proportional to $\Omega$ then the function $U(\Omega)=q\Omega^{-1}\in L_2(X, |\Omega|)$ is necessarily unbounded at some conical point $P_k$ and, therefore, does not belong to the domain of the Friedrichs extension of the operator $\Delta^{(0)}_{|\Omega|}$. Therefore,  the results of \cite{KokKorLeipzig} (where the explicit expression for the Friedrichs scalar Laplacian in the metric $|\Omega|$ was found) are not directly applicable. To include the holomorphic differentials into the domain of $\Delta^{(2)}_{|\Omega|}$ one has to consider some non standard self-adjoint extension of the scalar Laplacian.
 Self-adjoint extensions of the conical scalar Laplacian are described via asymptotics of the functions from their domains at the conical points. The extension of interest, the so called holomorphic extension $\Delta^{(0)}_{|\Omega|, \ hol}$ can be non-formally (the precise statement will be given later) described as follows. The functions $v$ from ${\cal D}(\Delta^{(0)}_{|\Omega|, \ hol})$ are subject to asymptotics
 $$v= \frac{H_k}{\zeta_k}+c_k+h_k\zeta_k+O(|\zeta_k|^{3/2})$$
 as $\zeta_k\to 0$, $k=1, \dots, 4g-4$, where
 \begin{equation}\label{dist}\zeta_k(P)=\left\{\int_{P_k}^P\sqrt{\Omega}\right\}^{2/3}\end{equation}
 is the so called distinguished local parameter near the conical point $P_k$.

 Using the results of L. Hillairet and the author  (\cite{HillKok1}, \cite{HillKok2}, see also \cite{KokLag}), we will establish a comparison formula (\ref{MAIN1}) for the determinants of the Friedrichs and the holomorphic extensions of the operator  $\Delta^{(0)}_{|\Omega|}$. This comparison formula slightly resembles heuristic formula (4.7) from \cite{Giddings} but is way more complicated. Formula (\ref{MAIN1}) together with results from \cite{KokKorLeipzig} imply an explicit expression for ${\rm det}\,  \Delta^{(0)}_{|\Omega|, \ hol}$.     Since the operator $\Delta^{(0)}_{\Omega, \ hol}$ is unitary equivalent to the self-adjoint  extension of the operator
$\Delta^{(2)}_{|\Omega|}$ whose domain contains all the holomorphic sections of the bundle $K^2$, we conclude that the usual $\zeta$-regularization of the determinant of the operator $\Delta^{(2)}_{|\Omega|}$  leads to the relation
\begin{equation}\label{DEF1}{\rm det}^{*}\,  \Delta^{(2)}_{|\Omega|}={\rm det}^{*}\,  \Delta^{(0)}_{|\Omega|, \ hol}\end{equation}
and get an explicit expression for the (modified, i. e. with zero modes excluded) determinant of the conical Laplacian in $K^2$. All this constitutes the first part of the present paper.

It is well-known that for smooth conformal metrics $\rho$ on $X$ the quantity \begin{equation}\label{Mumford0}\frac{{\rm det}\, \Delta^{(2)}_\rho}{{\rm det}||<q_i, q_j>_\rho||}\left(\frac{{\rm det}\,\Delta^{(0)}_\rho}{{\rm Area}\,(X, \rho){\rm det}\, \Im {\mathbb B}}\right)^{-13},\end{equation}
where ${\mathbb B}$ is the matrix of $b$-periods of $X$ and $||<q_i,  q_j>_\rho||$  is the Gram matrix for a basis $\{q_i\}$ of holomorphic quadratic differentials on $X$, is independent of the choice of the metric $\rho$ in the conformal class of $X$.
Moreover, the explicit expressions for this quantity (the modulus square of the so called Mumford measure on the moduli space) were found in \cite{Manin}, \cite{Sonoda}, \cite{Fay92}. Naive expectations that the same $\rho$-independence of (\ref{Mumford}) will hold if one extends the class of smooth conformal metrics to the class of conformal metrics with conical singularities fail: substituting the explicit expressions for
${\rm det}\Delta^{(2)}_{|\Omega|}$ and ${\rm det}\Delta^{(0)}_{|\Omega|}$ into (\ref{Mumford0}) gives a quantity which has nothing to do with any of the known explicit expressions for the Mumford measure; moreover the latter quantity seemingly  does not have any holomorphic factorization at all.

In \cite{EKZ} Eskin, Kontsevich and Zorich used an ingenious trick to define the determinant of the scalar Laplacian $\Delta^{(0)}_{|\omega|^2}$ in the conical metric given by the modulus square of a holomorphic one form $\omega$  on $X$ without any reference to the spectrum of the operator $\Delta^{(0)}_{|\omega|^2}$. Simple calculation based on the Burghelea-Friedlander-Kappeler gluing formula shows that their definition essentially (up to a constant) coincides with the usual $\zeta$-regularization of the (Friedrichs) scalar conical Laplacian. It turns out that Eskin-Kontsevich-Zorich (EKZ in what follows) regularization is also possible for the determinant of the operator $\Delta^{(2)}_{|\Omega|}$  although no relations between this regularization and the spectrum of the latter operator can be traced (to the best of our knowledge). In the second part of the present paper we find an explicit expression for the EKZ regularization of ${\rm det} \Delta^{(2)}_{|\Omega|}$ and show that substituting  this expression in (\ref{Mumford0}),  one gets a well-known formula for the Mumford measure.

{\bf Acknowledgements.} The author acknowledges numerous stimulating discussions with L. Hillairet and D. Korotkin. He especially thanks D. Korotkin for showing him the construction of the basis of quadratic differentials which is used in Section 3.1.   The research of the author was supported by NSERC. The preliminary version of the paper was written during the stay of the author in the Max Planck Institute for Mathematics in Bonn. The author thanks the Institute for hospitality and excellent working conditions.
\section{Holomorphic extension of scalar Laplacian in metric $|\Omega|$}
In this section we  derive an explicit expression for the determinant of the holomorphic extension of the scalar     Laplacian in metric $|\Omega|$ from the results of \cite{HillKok1}, \cite{HillKok2} and \cite{KokKorLeipzig}.

\subsection{Friedrichs and holomorphic extensions}
We start with description of the self-adjoint extensions of the symmetric operator $\Delta^{(0)}_{|\Omega|}$ initially defined on the set $C^\infty_0(X\setminus\{P_1, \dots, P_{4g-4}\})\subset L_2(X; |\Omega|)$.  Let ${\cal D}_{min}$ and  ${\cal D}_{max}$  be the domains of the closure of $\Delta^{(0)}_{|\Omega|}$ and the adjoint operator  $(\Delta^{(0)}_{|\Omega|})^*$ respectively. The following proposition (see, e. g., \cite{KokLag}, Appendix) gives an explicit description of the factor space ${\cal D}_{max}/{\cal D}_{min}$.
\begin{prop} In a vicinity of the point zero $P_k$ of the quadratic differential $\Omega$ functions $u$ from ${\cal D}_{max}$ are subject to the asymptotics
\begin{equation}\label{asy}
u(\zeta_k, \bar \zeta_k)=\frac{i}{\sqrt{2\pi}}L_k(u)\log|\zeta_k|+\frac{1}{4\pi}H_k(u)\frac{1}{\zeta_k}+\frac{1}{4\pi}A_k(u)\frac{1}{\bar \zeta_k}+\end{equation} $$\frac{i}{\sqrt{2\pi}}c_k(u)+\frac{1}{\sqrt{4\pi}}h_k(u)\zeta_k+\frac{1}{\sqrt{4\pi}}a_k(u)\bar\zeta_k+\chi(\zeta_k, \bar\zeta_k) v(\zeta_k, \bar \zeta_k)\,,
$$
where $\zeta_k$ is the distinguished local parameter (\ref{dist}) near $P_k$, $\chi$ is a cut-off function with support in a vicinity of $P_k$ and $v\in {\cal D}_{min}$. Any function $v$ from ${\cal D}_{min}$ is $O(|\zeta_k|^{3/2})$
near $P_k$.
\end{prop}
The next proposition can be used to describe all the self-adjoint extensions of the operator $\Delta^{(0)}_{|\Omega|}$.
\begin{prop} (see, e. g.,  \cite{HillKok2}, \cite{KokLag})
Let $u\in {\cal D}_{max}$ and let the row $$X(u)=(L_k(u), H_k(u), A_k(u), c_k(u), h_k(u), a_k(u))$$ consist of the coefficients in the asymptotics (\ref{asy}).
Choose representatives $u, v\in {\cal  D}_{max}$ of the classes $[u], [v]\in {\cal D}_{max}/{\cal D}_{min}$. Then
\begin{equation}\label{form}
\Omega([u], [v]):=\int_X \left( \left[(\Delta^{(0)}_{|\Omega|})^*u\right] \, v- u\,  \left[(\Delta^{(0)}_{|\Omega|})^*v\right]\right)|\Omega|=
\end{equation}
$$\sum_{k=1}^{4g-4}X_k(u)\left(\begin{matrix}0\ \ \ -I_3\\ I_3 \ \ \ \ 0   \end{matrix}    \right)X_k(v)^t\,.$$
\end{prop}
The s. a. extensions of the operator $\Delta^{0)}_{|\Omega|}$ are in one to one correspondence  with Lagrangian subspaces of ${\cal D}_{max}/{\cal D}_{min}$ with respect to symplectic form
$(u, v)\mapsto \Omega([u], [\bar v])$.
It is known the the Friedrichs extension $\Delta^{(0)}_{|\Omega|, F}$ corresponds to the Lagrangian subspace
$$L_k(u)=H_k(u)=A_k(u)=0\,.$$
Following L. Hillairet (\cite{Hillairet-French}), we define the {\it holomorphic extension}  $\Delta^{(0)}_{|\Omega|, hol}$ as the one corresponding to the lagrangian subspace
$$L_k(u)=A_k(u)=a_k(u)=0\,.$$
Thus, the function $u\in {\cal D}_{max}$ belongs to the domain of $\Delta^{(0)}_{|\Omega|, hol}$ if and only if the coefficients near the logarithm and the antiholomorphic terms $\bar \zeta_k^{-1}$ and $\bar \zeta_k$ in
(\ref{asy}) vanish for all $k=1, \dots, 4g-4$.
\begin{prop} The kernel of the Friedrichs extension $\Delta^{(0)}_{|\Omega|, F}$ is one-dimensional and consists of constants. The kernel of the holomorphic extension $\Delta^{(0)}_{|\Omega|, hol}$ has dimension $3g-3$ and consists of meromorphic functions with at most simple poles at the zeroes of the quadratic differential $\Omega$.
\end{prop}
{\bf Proof.} The first statement is trivial, let us prove the second. let $u\in {\rm Ker}\, \Delta^{(0)}_{|\Omega|, hol}$. Then
$$0=<\Delta^{(0)}_{|\Omega|, hol}u, u>=\lim_{\epsilon\to 0}\left[\iint_{X_\epsilon}\bar u \partial\bar \partial u\right]=
-\iint_X|\bar \partial u|^2+\lim_{\epsilon\to 0}\int_{\partial X_{\epsilon}}\bar u\bar\partial u,$$
where $X_\epsilon$ is obtained from $X$ cutting out disks of radius $\epsilon$ (in the distinguished local parameters) around all the points $P_k$.
Using asymptotics $u\sim H_k/\zeta_k+c_k+h_k\zeta_k+O(|\zeta_k|^{3/2})$ at $P_k$, one easily shows that the limit of the last contour integral is equal to zero and, therefore, $u$ is holomorphic in $X\setminus \{P_1, \dots, P_{4g-4}\}$.
It is clear from the same asymptotics  that the poles of $u$ at $P_k$ are at most simple.  Therefore, the quadratic differential $u\Omega$ is holomorphic. On the other hand any function of the form $q/\Omega$, where $q$ is a holomorphic quadratic differential, belongs to  $u\in {\rm Ker}\, \Delta^{(0)}_{|\Omega|, hol}$. Thus, the dimension of the kernel coincides with the dimension of the space of holomorphic quadratic differentials and is equal to $3g-3$. $\square$.

\subsection{$S$-matrix of the conical surface $X$} Here we introduce the so-called $S$-matrix of a conical surface $X$ which enters as the main ingredient into the formula relating the zeta-regularized determinants of two different self-adjoint extensions of the scalar Laplacian on $X$ (\cite{HillKok1}).
Let $\lambda$ be the spectral parameter; assuming that $\lambda$ do not belong the (discrete, nonnegative) spectrum of the Friedrichs extension $\Delta^{(0)}_{|\Omega|, F}$, define {\it the special growing solutions}
$G_{1/\zeta_k}(\, \cdot\,; \lambda)$, $G_{1/\bar\zeta_k}(\, \cdot\,; \lambda)$, $G_{\log|\zeta_k|}(\, \cdot\,; \lambda)$
of the equation
$$(\Delta^{(0)}_{|\Omega|})^*G-\lambda G=0$$
specifying their asymptotics near the the conical points $P_k$ of the surface $X$:
$$G_{1/\zeta_k}(P; \lambda)=\frac{1}{\zeta_k(P)}+O(1)\ \ \ \ \text{as}\ \ P\to P_k,$$
$$G_{1/\zeta_k}(P; \lambda)=O(1)\ \ \ \ \text{as}\ \ P\to P_l;\ \  l\neq k;$$

$$G_{1/\bar \zeta_k}(P; \lambda)=\frac{1}{\overline{\zeta_k(P)}}+O(1)\ \ \ \ \text{as}\ \ P\to P_k,$$
$$G_{1/\bar \zeta_k}(P; \lambda)=O(1)\ \ \ \ \text{as}\ \ P\to P_l;\ \  l\neq k;$$

$$G_{\log|\zeta_k|}(P; \lambda)=\log|\zeta_k(P)|+O(1)\ \ \ \ \text{as}\ \ P\to P_k,$$
$$G_{\log|\zeta_k|}(P; \lambda)=O(1)\ \ \ \ \text{as}\ \ P\to P_l; \ \ l\neq k.$$
The entries $S^{\alpha \beta}(\lambda)$
$$\alpha=\log|\zeta_k|, \ \ 1/\zeta_k,\ \  1/\bar \zeta_k,\ \   k=1, \dots, 4g-4;$$ $$\beta=1_l,\ \  \zeta_l, \ \ \bar \zeta_l,\ \  l=1, \dots, 4g-4$$
 of the $S$-matrix of the conical surface $(X, |\Omega|)$ are defined via the asymtotics
$$ G_{1/ \zeta_k}(P; \lambda)=\delta_{kl}\frac{1}{\zeta_l(P)}+S^{1/\zeta_k, 1_l}(\lambda)+S^{1/\zeta_k, \zeta_l}(\lambda)\zeta_l(P)+S^{1/\zeta_k, \bar\zeta_l}(\lambda)\bar\zeta_l(P)+O(|\zeta_l(P)|^{3/2}), $$
$$ G_{1/   \bar\zeta_k}(P; \lambda)=\delta_{kl}\frac{1}{\overline{\zeta_l(P)}}+S^{1/\bar\zeta_k, 1_l}(\lambda)+S^{1/\bar\zeta_k, \zeta_l}(\lambda)\zeta_l(P)+S^{1/\bar\zeta_k, \bar\zeta_l}(\lambda)\bar\zeta_l(P)+O(|\zeta_l(P)|^{3/2}), $$
$$ G_{\log|\zeta_k|}(P; \lambda)=\delta_{kl}\log|\zeta_l(P)|+S^{\log|\zeta_k|, 1_l}(\lambda)+S^{\log|\zeta_k|, \zeta_l}(\lambda)\zeta_l(P)+$$$$S^{\log|\zeta_k|, \bar\zeta_l}(\lambda)\bar\zeta_l(P)+O(|\zeta_l(P)|^{3/2}), $$
as $P\to P_l$.
We refer to \cite{HillKok1}, \cite{HillKok2}, \cite{HKK-Trans} for the study of analytic properties properties of the entries of $S(\lambda)$ (see also \cite{KokLag} for an improved exposition).
In what follows we will only deal with the $(4g-4)\times(4g-4)$ block $T(\lambda)$ of the matrix $S(\lambda)$:
\begin{equation}\label{defT}T(\lambda)=||S^{1/\zeta_k, \bar \zeta_l}(\lambda)||_{k, l=1, \dots, 4g-4}\,.\end{equation}
The matrix $T(\lambda)$ is known to be analytic in a vicinity of $\lambda=0$, moreover, one has the relation
\begin{equation}\label{S-explicit}
T(0)=-\pi||B(P_k, P_l)||_{k, l=1, \dots, 4g-4},
\end{equation}
where
$$B(R, Q)=\sum_{i, j=1}^g(\Im {\mathbb B})_{ij}^{-1}v_i(R)\overline{v_l(Q)}$$
is the Bergman kernel of the Riemann surface $X$. Here ${\mathbb B}$ is the matrix of $b$-periods of $X$, $\{v_i\}_{i=1, \dots, g}$ is the basis of normalized holomorphic differentials on $X$ and $v_j(P_l)$ means the value
of $v_j$ at $P_l$ computed in the distinguished local parameter $\zeta_l$:
$$v_j=f(\zeta_l)d\zeta_l\ \ \ \text{near}\ \ P_l;  \ \ \ \ \  v_j(P_l):=f(0)\,.$$
The following proposition almost immediately follows from Proposition 2 from \cite{HillKok2}.
\begin{prop} The matrix $T(0)$ has rank $g$.
\end{prop}
{\bf Proof.} The matrix $T(0)$ is the Gram matrix of the $4g-4$ vectors
$${\bf V}_l=\left(v_1(P_l), \dots, v_g(P_l)\right); \ \ l=1, \dots, 4g-4$$
from ${\mathbb C}^g$ with respect to the Hermitian product
$$<{\bf V}, {\bf W}>=\sum_{i, j}^g(\Im {\mathbb B})_{ij}^{-1}V_i\bar W_j$$
and, therefore, ${\rm rank}\, T(0)\leq g$. Assume that  ${\rm rank}\, T(0)< g$. Then Proposition 2 from \cite{HillKok2} implies that for any subset $\{P_{i_1}, \dots, P_{i_{2g-2}}\}$ of the set of zeroes $P_1, \dots, P_{4g-4}$ of the quadratic differential $\Omega$ the divisor $P_{i_1}+\cdots P_{i_{2g-2}}$ belongs to the canonical class. Considering two subsets $\{P_1, P_2, \dots, P_{2g-3}, P_{2g-2}\}$ and
$\{P_1, P_2, \dots, P_{2g-3}, P_{2g-1}\}$, one gets existence of a meromorphic function with divisor $P_{2g-2}-P_{2g-1}$ which is impossible. $\square$
\subsection{Comparison formula for determinants of Friedrichs and holomorphic extensions}
The following Proposition is a version of Theorem 2 from \cite{HillKok1} and Proposition 3 from \cite{HillKok2}, its proof coincides with the one from \cite{HillKok2} almost verbatim. The absolute constant $C_g$ below admits explicit calculation (cf. Propositions 7 and 8 from \cite{KokLag}) but this is of no interest for the purpose of the present paper.

\begin{prop} Let the spectral parameter $\lambda$ do not belong to the union of spectra of the Friedrichs, $\Delta^{(0)}_{|\Omega|, F}$,  and the holomorphic, $\Delta^{(0)}_{|\Omega|, hol}$, extensions of the scalar laplacian.
Then the $\zeta$-regularized determinants of the operators  $\Delta^{(0)}_{|\Omega|, F}-\lambda$  and $\Delta^{(0)}_{|\Omega|, hol}- \lambda$ are related as follows:
\begin{equation}\label{comparison1}
{\rm det}\, (\Delta^{(0)}_{|\Omega|, hol}- \lambda)=C_g{\rm det}\, T(\lambda){\rm det}\, (\Delta^{(0)}_{|\Omega|, F}- \lambda)\,
\end{equation}
where $T(\lambda)$ is from (\ref{defT}) and the constant $C_g$ depends only on the genus of $X$.
\end{prop}

Since ${\rm dim}\,{\rm Ker}\,\Delta^{(0)}_{|\Omega|, F}=1$ and ${\rm dim}\,{\rm Ker}\,\Delta^{(0)}_{|\Omega|, hol}=3g-3$, one gets the following relation between the modified (i. e. with zero modes excluded) determinants of the two extensions ${\rm dim}\,{\rm Ker}\,\Delta^{(0)}_{|\Omega|, F}$ and ${\rm dim}\,{\rm Ker}\,\Delta^{(0)}_{|\Omega|, hol}$.
\begin{thm}\label{T1} One has the relation
\begin{equation}\label{MAIN1}{\rm det}^*\, \Delta^{(0)}_{|\Omega|, hol}=\frac{C_g}{(3g-4)!}\left(\frac{d}{d\lambda}\right)^{3g-4}{\rm det}\, T(\lambda)\Big|_{\lambda=0}{\rm det}^*\, \Delta^{(0)}_{|\Omega|, F}\,.\end{equation}
\end{thm}
Denoting by ${\cal N}_{3g-4}$ the set of $(4g-4)$-tuples ${\bf n}=(n_1, \dots, n_{4g-4})$ of non-negative integers such that $n_1+\dots +n_{4g-4}=3g-4$, one has
\begin{equation}\label{diffDet}
\left(\frac{d}{d\lambda}\right)^{3g-4}{\rm det}\, T(\lambda)\Big|_{\lambda=0}=
\end{equation}
$$\sum_{{\bf n}\in {\cal N}_{3g-4}} \frac{(3g-4)!}{n_1!\,n_2!\,\dots\, n_{4g-4}!}\left|\begin{matrix}\left(\frac{d}{d\lambda}\right)^{n_1}S^{1/\zeta_1, \bar \zeta_1}(0) \ \ \  \dots \ \ \ \left(\frac{d}{d\lambda}\right)^{n_{4g-4}}S^{1/\zeta_{1}, \bar \zeta_{4g-4}}(0)  \\ \dots \\\left(\frac{d}{d\lambda}\right)^{n_1}S^{1/\zeta_{4g-4}, \bar \zeta_1}(0) \dots \left(\frac{d}{d\lambda}\right)^{n_{4g-4}}S^{1/\zeta_{4g-4}, \bar \zeta_{4g-4}}(0)
   \end{matrix}\right|   $$

All the entries of the determinants from the right hand side of (\ref{diffDet}) can be explicitly computed (see \cite{HillKok1}, formulae (4.12), (4.9)  and (2.8); and \cite{KokLag}, formula (2.20)). Namely,
introduce the meromorphic differential on $X$ of the third kind
$$\Omega_{p-q}(z)=\int_p^qW(z,\, \cdot\, )-2\pi \sqrt{-1}\sum_{i, j=1}^g(\Im{\mathbb B})^{-1}_{ij}v_\alpha(z)\int_p^qv_\beta\,$$
with simple poles at $p$ and $q$ with residues $1$ and $-1$ and purely imaginary periods; here $W(\, \cdot\,,\, \cdot\,)$ is the canonical meromorphic bidifferntial on $X$.
Then (see \cite{KokLag}, Proposition 3) the special growing solution $G_{1/\zeta_k}(P, \lambda)$ admits holomorphic continuation to $\lambda=0$ and one has
\begin{equation}
G_{1/\zeta_k}(y; 0)=-\frac{1}{\int_X|\Omega|}\int_X\Omega_{y-\cdot}(P_k)|\Omega(\cdot)|,
\end{equation}
where the value of the $1$-differential at $P_k$ is computed in the distinguished local parameter;
moreover, the function $G_{1/\zeta_k}(y; 0)$ is orthogonal to $1$ in $L_2(X, |\Omega|)$.
Thus, formulae (4.11) and (4.12) from \cite{HillKok2} imply
the relation
\begin{equation}\label{diffentry}
\left(\frac{d}{d\lambda}\right)^{n}S^{1/\zeta_l, \bar \zeta_k}(0)=\int_X \left[\left(\Delta^{(0)}_{|\Omega|, F}\Big|_{1^\bot}\right)^{1-n}G_{1/\zeta_l}(\,\cdot\,, 0)\right]\overline{ G_{1/\zeta_k}(\,\cdot\,, 0)}|\Omega(\,\cdot\,)|\,
\end{equation}
where the inverse operator
$\left(\Delta^{(0)}_{|\Omega|, F}\Big|_{1^\bot}\right)^{-1}:1^\bot \to 1^\bot$
acts as
$u\mapsto v$, where
$$v(x)=\int_X{\frak{G}}(x, \,\cdot\,)u(\,\cdot\,)|\Omega(\,\cdot\,)|$$
and
$${\frak{G}}(x, y)=\frac{1}{2\pi \left(\int_X|\Omega|\right)^2 }\int_{q\in X}|\Omega(q)|\int_{p\in X}|\Omega(p)|\Re\int_p^x\Omega_{y-q}(\,\cdot\,)$$
is the Green function of the Friedrichs Laplacian (see \cite{KokLag}, formula (2.5)).

\begin{Remark}Following the discussion from the Introduction, one can consider
Theorem \ref{T1} together with formulae (\ref{diffDet}) and (\ref{diffentry}) as a rigorous counterpart of the heuristic formula (4.7) from \cite{Giddings}.
\end{Remark}

\subsection{Explicit formula for the determinant of the Friedrichs extention}
Here we recall an explicit formula for ${\rm det}^*\, \Delta^{(0)}_{|\Omega|, F}$ obtained in \cite{KokKorLeipzig}; the latter formula together with Theorem \ref{T1} and \ref{DEF1} implies an explicit expression for the $\zeta$-regularized determinant of the Laplacian in $K^2$ corresponding to the flat conical metric $|\Omega|$.

Let $Q_g(1, \dots, 1)$  ($4g-4$ units) be the the moduli space of  pairs $(X, \Omega)$, where $X$ is a compact Riemann surface of genus $g$ and $\Omega$  is a quadratic holomorphic differential on $X$ with ($4g-4$) simple zeroes. Following the improved presentation from more recent paper \cite{KoroZograf},  define the Bergman tau-function $\tau(X, \Omega)$ on $Q_g(1, \dots, 1)$ via the relation
\begin{equation}\label{tau}
\tau^{24}(X, \Omega)=\frac{\left( \left[\sum_{i=1}^gv_j(P)\frac{\partial}{\partial w_j}\right]^g\Theta(w; {\mathbb B})\Big|_{w=K^P} \right)^{16} \prod_{k<l}E(P_k, P_l)}    {e^{4\pi\sqrt{-1}<{\mathbb B}Z+4K^P, Z>}{\cal W}^{16}(P)\prod_kE(P, P_k)^{4(g-1)}}\,
\end{equation}
where
\begin{itemize}
\item ${\mathbb B}$ is the matrix of $b$-periods of $X$; $\Theta(\,\cdot\,, {\mathbb B}$ is the Riemann theta-function,
\item $P_k$, $k=1, \dots, 4g-4$ are the zeroes of the quadratic differential $\Omega$.
\item $P$ is an arbitrary base point on $X$ (in fact, expression (\ref{tau}) is     $P$-independent),
\item ${\cal W}$ is the Wronskian of basic holomorphic differentials $v_1, \dots, v_g$,
\item $E(\cdot, \cdot)$ is the prime form,
\item $K^P$ is the vector of Riemann constants,
\item the vector $Z$ from $\left(\frac{1}{2}{\mathbb Z}\right)^g$  is defined via
$$\frac{1}{2}{\cal A}_P((\Omega))+2K^{P}={\mathbb B}Z+Z'\,$$
where ${\cal A}_P$ is the Abel map with base $P$ and $(\Omega)=P_1+\dots+P_{4g-4}$ is the divisor of $q$; $Z'\in \left(\frac{1}{2}{\mathbb Z}\right)^g$,
\item all the values of tensor (the basic holomorphic differentials)  and tensor-like (the prime-form) objects at $P_k$ are computed in the distinguished local parameter
$$\zeta_k(Q)=\left\{\int_{P_k}^Q\sqrt{\Omega}\right\}^{2/3}$$
near $P_k$,
\item all the values of tensor (the Wronskian)  and tensor-like (the prime-form) objects at $P$ are computed in the same holomorphic local parameter at $P$; since expression (\ref{tau}) is a scalar w. r. t. $P$, the result does not depend on the choice of this local parameter.
\end{itemize}
\begin{thm}\label{Leipzig} Let $(X, \Omega)$ be an element of the space $Q_g(1, \dots, 1)$.  The $\zeta$-regularized (modified) determinant of the Friedrichs extension of the scalar Laplacian in the metric $|\Omega|$ on $X$ can be expressed through
the Bergman tau-function $\tau$ on $Q_g(1, \dots, 1)$  from (\ref{tau}) as follows:
\begin{equation}\label{dettau}
{\rm det}^*\,\Delta^{(0)}_{|\Omega|, F}=C(g)\left(\int_X|\Omega|\right)\, {\rm det}\Im{\mathbb B}\, |\tau|^2\,,\end{equation}
where the constant $C(g)$ depends only on the genus $g$.
\end{thm}

\subsection{EKZ-regularization of the scalar Laplacian and BFK gluing formulae}
Here we show that the regularization  of the determinant of the Friedrichs scalar Laplacian $\Delta^{(0)}_{|\Omega|, F}$ on $X$ proposed by Eskin, Kontsevich and Zorich in \cite{EKZ} essentially coincides with usual $\zeta$-regularization. This fact is derived using the Bourghlea-Fiedlander-Kappeler (BFK in what follows) gluing formalae. First, we recall the EKZ construction. It should be noted that in the original paper \cite{EKZ} this construction was introduced to regularize the determinant of the Laplacian $\Delta^{(0)}_{|v|^2, F}$ considered as a function on the moduli space $H_g(k_1, \dots, k_m)$ of pairs $(X, v)$, where $X$ is a compact Riemann surface and $v$ is a holomorphic one form with $M$ zeroes of multiplicities $k_1, \dots, k_M$; here we apply this construction without any changes to the pairs $(X, \Omega) \in Q_g(1, \dots, 1)$.

Choose a pair $(X_0, \Omega_0)$ from the space $Q_g(1, \dots, 1)$ a let $(X, \Omega)$ be any other element of $Q_g(1, \dots, 1)$. In a vicinity of any zero $P_k$ of the quadratic differential $\Omega$ the metric $|\Omega|$ coincides with the metric of the standard round cone of the angle $3\pi$,
$$|\Omega|=\frac{9}{4}|\zeta_k||d\zeta_k|^2$$
(one can certainly get rid  of the factor $9/4$ just including the corresponding factor in the definition of the distinguished local parameter $\zeta_k$  (\ref{dist})).
Smoothing all these standard round cones in $\epsilon$-vicinities of their tips $P_1, \dots, P_{4g-4}$ one gets a {\it smooth conformal} metric $|\Omega|^{(\epsilon)}$ on $X$. One can perform the smoothing in such  a way that the total area of $X$ does not change, i . e. $\int_X|\Omega|=\int_X|\Omega|^{(\epsilon)}$. The following key observation belongs to Eskin, Kontsevich and Zorich.
\begin{lem}\label{EKZlemma}
The quantity
\begin{equation} \label{ratio} R^{(0)}_\epsilon(X, \Omega):=\frac{{\rm det}^*\Delta^{(0)}_{|\Omega|^{(\epsilon)}}}{{\rm det}^*\Delta^{(0)}_{|\Omega_0|^{(\epsilon)}}}\end{equation}
is independent of $\epsilon$ for sufficiently small positive $\epsilon$. Here both the determinants of the Laplacians in smooth metrics $|\Omega_0|^{(\epsilon)}$ (on $X_0$) and $|\Omega|^{(\epsilon)}$ (on $X$)  are defined via usual $\zeta$-regularization.
\end{lem}
This lemma easily follows from the Polyakov formula relating the determinants of the Laplacians in two smooth conformally equivalent metrics on a given Riemann surface. In the next section we will prove its complete analog
for the Laplacians acting in $K^2$ giving all the needed details.
\begin{defi} Let $(\Omega, X)\in Q_g(1, \dots, 1)$. Define the EKZ-regularization of the determinant of the operator $\Delta^{(0)}_{|\Omega|}$ as
\begin{equation}\label{EKZ-def1}
{\rm det}_{EKZ} \Delta^{(0)}_{|\Omega|}:=\lim_{\epsilon \to 0}R^{(0)}_\epsilon(X, \Omega)\,.
\end{equation}
\end{defi}
\begin{Remark}We remind the reader that the Laplacians in smooth metrics are essentially self-adjoint, so no choice of the self-adjoint extension is required in the definition of $R^{(0)}_\epsilon(X, \Omega)$ as well as
in that of ${\rm det}_{EKZ} \Delta^{(0)}_{|\Omega|}$; in a sense the latter quantity is defined without any use of the spectrum of the operator $\Delta^{(0)}_{|\Omega|}$.
\end{Remark}

The following Proposition shows that for the scalar conical Laplacians the EKZ regularization of the determinant essentially coincides with the usual $\zeta$-regularization.
\begin{prop}\label{ravny} One has the relation
\begin{equation}\label{identify}
{\rm det}_{EKZ} \Delta^{(0)}_{|\Omega|}=C {\rm det}^*\Delta^{(0)}_{|\Omega|, F}\,
\end{equation}
where the constant $C$ is independent of the point $(X, \Omega)\in Q_g(1, \dots, 1)$.
\end{prop}
{\bf Proof.} The following argument is very similar to the proof of Proposition 1 from \cite{KokProcAMS} (and we refer to this reference for all the details omitted): one uses the fact that in vicinities ($D_k(\epsilon)$) of conical points  all the Riemannian manifolds $(X, |\Omega|)$ are isometric to the vicinity  $|x|<\epsilon$ of the tip of the standard round cone $({\mathbb C}, |x||dx|^2)$; after smoothing all these vicinities turn into the standard smooth "cap" $(\{|x|\leq \epsilon\}, \rho(x, \bar x)|dx|^2)$ with smooth positive $\rho$. One can assume that the metric $\rho(x, \bar x)|dx|^2$ coincides with $|x||dx|^2$ for $0<\epsilon_1<|zx|\leq \epsilon$ and, therefore the lengths of the boundary circles $|x|=\epsilon$ do not change after smoothing. We also remind the reader that the metrics $|\Omega|$ and $|\Omega|^{(\epsilon)}$ have the same volume.

Let $\Sigma=X\setminus \cup_{k=1}^{4g-4}D_k(\epsilon)$, $\Gamma=\partial \Sigma$. Denote by $l(\Gamma)$  the length of $\partial \Gamma$ (in the metric $|\Omega|$). Let also ${\cal N}$ be the Neumann jump operator on the contour $\Gamma$ for $\Delta^{(0)}_{|\Omega|^{(\epsilon)}}$.
Let $\Sigma_0$, $\Gamma_0$, ${\cal N}_0$ be the same objects defined for the reference surface $(X_0, |\Omega_0|^{(\epsilon)})$.
In the sequel we denote by $(\Delta|{\cal U})$ the operator of the Dirichlet boundary value problem for an operator $\Delta$ in a domain ${\cal U}$.

Applying the BFK gluing formulae
\begin{equation}\label{BFK1}
{\rm det}^*\Delta^{(0)}_{|\Omega|^\epsilon}=\left\{\prod_{k=1}^{4g-4}{\rm det}\,(\Delta^{(0)}_{|\Omega|^{(\epsilon)}}; D_k(\epsilon))\right\}{\rm det}(\Delta^{(0)}_{|\Omega|}|\Sigma){\rm det}^*{\cal N} {\rm Area}(X, |\Omega|) l(\Gamma)^{-1}
\end{equation}
to the numerator and denominator of (\ref{ratio}),
one gets
\begin{equation}\label{auxil}R^{(0)}_\epsilon(X, \Omega)=\frac{{\rm det}(\Delta^{(0)}_{|\Omega|}|\Sigma){\rm det}^*{\cal N} {\rm Area}(X, |\Omega|) \l(\Gamma)^{-1}}
{{\rm det}(\Delta^{(0)}_{|\Omega_0|}|\Sigma_0){\rm det}^*{\cal N}_0 {\rm Area}(X_0, |\Omega_0|)\l(\Gamma_0)^{-1}}\,.\end{equation}
Multiplying the right hand side of (\ref{auxil}) by
$$1=\frac{\left\{\prod_{k=1}^{4g-4}{\rm det}\,(\Delta^{(0)}_{|\Omega|}| D_k(\epsilon))\right\}}{\left\{\prod_{k=1}^{4g-4}{\rm det}\,(\Delta^{(0)}_{|\Omega_0|}| D_k(\epsilon))\right\}}$$
observing that the determinant of the Neumann jump operator does not change if one passes from the metric $|\Omega|^{(\epsilon)}$ to the metric $|\Omega|$ and then making use of the following BFK formula
for the operator $\Delta^{(0)}_{|\Omega|, F}$
\begin{equation}
{\rm det}^*\Delta^{(0)}_{|\Omega|, F}= \left\{\prod_{k=1}^{4g-4}{\rm det}\,(\Delta^{(0)}_{|\Omega|}; D_k(\epsilon))\right\}{\rm det}(\Delta^{(0)}_{|\Omega|}|\Sigma){\rm det}^*{\cal N} ({\rm Area}(X, |\Omega|) l(\Gamma)^{-1}\,,
\end{equation}
one gets the relation
$$R^{(0)}_\epsilon(X, \Omega)=\frac{{\rm det}^*\Delta^{(0)}_{|\Omega|, F}}{{\rm det}^*\Delta^{(0)}_{|\Omega_0|, F}}$$
and, therefore, the statement of the proposition. $\square$

\section{EKZ-regularization of the determinant of $\Delta^{(2)}$}

\subsection{Canonical covering and Korotkin-Zograf basis}
To a pair $(X, \Omega)$ from $Q_g(1, \dots, 1)$ one associates a ramified two-fold covering $\pi: \hat{X}\to X$, where $\hat{X}$ (the so-called canonical cover) is a compact Riemann surface of genus $4g-3$. The quadratic differential $Q$ being lifted to $\hat{X}$ becomes the square of a holomorphic differential $v$ on $\hat {X}$. The branch points of the covering coincide with $4g-4$ zeroes of the quadratic differential $Q$; the differential $v$ has double zeroes at preimages of the branch points. Let $\mu: \hat{X}\to \hat{X}$ be the holomorphic involution interchanging the sheets of the covering. Following \cite{Fay1} introduce a canonicle basis of $a$ and $b$-cycles on $\hat{X}$:
$$\{A_l\}_{l=1}^{4g-3}= \{a_\alpha\}_{\alpha=1}^g\cup \{\mu(a_\alpha)\}_{\alpha=1}^g\cup \{\tilde a_k\}_{k=1}^{2g-3}$$
$$\{B_l\}_{l=1}^{4g-3}= \{b_\alpha\}_{\alpha=1}^g\cup \{\mu(b_\alpha)\}_{\alpha=1}^g\cup \{\tilde b_k\}_{k=1}^{2g-3}$$
where the cycles $(a_\alpha, \mu(a_\alpha), b_\alpha, \mu(b_\alpha))$ are obtained via lifting of a canonical basis of $a$ and $b$-cycles on the base $X$ of the covering and the cycles $\tilde a_k$, $\tilde b_k$ satisfy the relations
$$\mu(\tilde a_k)+\tilde a_k=\mu(\tilde b_k)+\tilde b_k=0$$
in $H^1(X, {\mathbb Z})$.
Let
$$\{V_l\}_{l=1}^{4g-3}=\{v_\alpha\}_{\alpha=1}^g \cup\{ v^{(\mu)}_\alpha\}_{\alpha=1}^g\cup \{w_k\}_{k=1}^{2g-3}$$
 be the corresponding basis of normalized (i. e. $\int_{A_m}V_n=\delta_{mn}$) holomorphic differentials on $\hat{X}$:
Following Korotkin and Zograf (\cite{KoroZograf}, see also \cite{BKN-invent}), introduce $3g-3$ holomorphic differentials $W_1, \dots, W_{3g-3}$  on $\hat{X}$ via
$$W_\alpha=v_\alpha-v^{(\mu)}_\alpha; \alpha =1, \dots, g$$
$$W_{g+k}=w_k; k=1, \dots, 2g-3\,.$$
Clearly, one has
$$\mu^*W_k=-W_k; k=1, \dots, 3g-3\,.$$
On the other hand the differential $v=\sqrt{\pi^*\Omega}$ satisfies
$$\mu^{*}v=-v$$ and, therefore,
the quadratic differentials $vW_k; k=1, \dots, 3g-3$ are invariant with respect to the holomorphic involution $\mu$ on $\hat{X}$, and, therefore, give rise to quadratic differentials $q_1, \dots, q_{3g-3}$ on $X$.

\begin{defi}The basis $q_1, \dots, q_{3g-3}$ in $H^0(X, K^2)$ is called Korotkin-Zograf basis of holomorphic quadratic differentials on $X$.
\end{defi}

\subsection{EKZ regularization}
Let $(X, \Omega)$ be a point of the stratum $Q_g(1, \dots, 1)$ and let $\epsilon>0$ be sufficiently small.
 As in Section 2 introduce the {smooth} conformal metric $|\Omega|^{(\epsilon)}$ on $X$ which coincides with $|\Omega|$ outside the disks $|x_m|\leq \epsilon$, where
$x_m(P)=\left(\int_{P_m}^P\sqrt{\Omega}\right)^{2/3}$ is the distinguished local parameter in a vicinity of $P_m$. One can assume that $\int_X|\Omega|^{(\epsilon)}=\int_X|\Omega|$.

Denote by ${\cal G}({q_\alpha}; |\Omega|^{(\epsilon)}; X, Q)$ the Gram determinant of the scalar products w. r. t. the metric $|\Omega|^{(\epsilon)}$:
$${\cal G}({q_\alpha}; |\Omega|^{(\epsilon)}; X, Q)={\rm det}\Big|\Big|\int_X \frac{q_\alpha q_\beta}{|\Omega|^{(\epsilon)}}\Big|\Big|$$
(here $|\Omega|^{(\epsilon)}$ is considered as a $(1, 1)$-form).

Choose a reference point $(X_0, \Omega_0)$ of the stratum $Q_g(1, \dots, 1)$.

Let
$$R^{(2)}_\epsilon(X, \Omega)=\frac{{\rm det}\Delta^{(2)}_{|\Omega|^{(\epsilon)}}(X, \Omega)\left({\cal G}({q_\alpha}; |\Omega|^{(\epsilon)}; X, Q)\right)^{-1}                        }{{\rm det}\Delta^{(2)}_{|\Omega_0|^{(\epsilon)}}(X_0, \Omega_0)\left({\cal G}({q_\alpha}; |\Omega_0|^{(\epsilon)}; X_0, \Omega_0)\right)^{-1} }\,.$$

The following statement is an analog of Lemma \ref{EKZlemma}.

\begin{lem}\label{EKZ}  For sufficiently small $\epsilon$ the quantity $R^{(2)}_\epsilon(X, \Omega)$ is $\epsilon$-independent.
\end{lem}
{\bf Proof.} The statement of the Lemma almost immediately follows from Proposition 3.8 from
\cite{Fay92}.

Let $\rho_{\epsilon_k}=\left\{|\Omega|^{(\epsilon_{k})}\right\}^{-1/2}$; $h_{\epsilon_{k}}=\left\{|\Omega|^{(\epsilon_{k})}\right\}^{-2}$
 and
 $\rho_{\epsilon_{k}}^0=\left\{|\Omega_0|^{(\epsilon_{k})}\right\}^{-1/2}$; $h_{\epsilon_{k}}^0=\left\{|\Omega_0|^{(\epsilon_{k})}\right\}^{-2}$; $k=1, 2$.

According to Proposition 3. 8 from \cite{Fay92} one has
$$\log R_{\epsilon_2}^{(2)}-\log R_{\epsilon_1}^{(2)}=
\log \frac{{\rm det}\Delta^{(2)}_{|\Omega|^{(\epsilon_2)}}(X, \Omega)\left({\cal G}({q_\alpha}; |\Omega|^{(\epsilon_2)}; X, Q)\right)^{-1}  }    {{\rm det}\Delta^{(2)}_{|\Omega|^{(\epsilon_1)}}(X, \Omega)\left({\cal G}({q_\alpha}; |\Omega|^{(\epsilon_1)}; X, Q)\right)^{-1}  }
\ - $$ $$\log\frac{{\rm det}\Delta^{(2)}_{|\Omega_0|^{(\epsilon_2)}}(X_0, \Omega_0)\left({\cal G}({q_\alpha}; |\Omega_0|^{(\epsilon_2)}; X_0, \Omega_0)\right)^{-1}  }  {  {\rm det}\Delta^{(2)}_{|\Omega_0|^{(\epsilon_1)}}(X_0, \Omega_0)\left({\cal G}({q_\alpha}; |\Omega_0|^{(\epsilon_1)}; X_0, \Omega_0)\right)^{-1}}=
$$
$$
\frac{1}{\pi}\int_X\left\{(\frac{1}{3}\log\frac{\rho_{\epsilon_2}}{\rho_{\epsilon_1}}-\frac{1}{2}\log \frac{h_{\epsilon_2}}{h_{\epsilon_1}})\partial^2_{z\bar z}\log \rho_{\epsilon_1}\rho_{\epsilon_2}+
(-\frac{1}{2}\log \frac{\rho_{\epsilon_2}}{\rho_{\epsilon_1}}+\frac{1}{2}\log \frac{h_{\epsilon_2}}{h_{\epsilon_1}})\partial^2_{z\bar z}\log\frac{h_{\epsilon_2}}{h_{\epsilon_1}}
\right\}\widehat{dz}\ -
$$
$$
\frac{1}{\pi}\int_X\left\{(\frac{1}{3}\log\frac{\rho^0_{\epsilon_2}}{\rho^0_{\epsilon_1}}-\frac{1}{2}\log \frac{h^0_{\epsilon_2}}{h^0_{\epsilon_1}})\partial^2_{z\bar z}\log \rho^0_{\epsilon_1}\rho^0_{\epsilon_2}+
(-\frac{1}{2}\log \frac{\rho^0_{\epsilon_2}}{\rho^0_{\epsilon_1}}+\frac{1}{2}\log \frac{h^0_{\epsilon_2}}{h^0_{\epsilon_1}})\partial^2_{z\bar z}\log\frac{h^0_{\epsilon_2}}{h^0_{\epsilon_1}}
\right\}\widehat{dz}\,.
$$
Assuming $\epsilon_2>\epsilon_1$, one has  $|\Omega|^{(\epsilon_1)}=|\Omega|^{(\epsilon_2)}$ and
$|\Omega_0|^{(\epsilon_1)}=|\Omega_0|^{(\epsilon_2)}$ outside the $\epsilon_2$-disks around the conical points. So the integration in the r. h. side of the above formula goes only these $\epsilon_2$-disks. Inside these disks the expressions under the two integrals in the r. h. s. are the same (the parts of the flat surfaces $(X, \Omega)$ and $(X_0, \Omega_0)$ inside these disks are isometric) and, therefore,
$\log R_{\epsilon_2}^{(2)}-\log R_{\epsilon_1}^{(2)}=0$. $\square$

\begin{defi}

Define the EKZ-regularization of the determinant of the Laplacian acting on quadratic differentials in the singular metric $|\Omega|$ as
\begin{equation}\label{detreg2}{\rm det}_{EKZ}\Delta^{(2)}(X, \Omega):=\left(\lim_{\epsilon\to 0}R^{(2)}_\epsilon (X, \Omega)\right) {\rm det} \Big| \Big| \int_X\frac{q_\alpha \bar q_\beta} {|\Omega|}    \Big| \Big| \,. \end{equation}
\end{defi}

\begin{Remark} It can be shown that the determinant of the matrix
$$\Big| \Big| \int_X\frac{q_\alpha \bar q_\beta} {|\Omega|}    \Big| \Big|=\Big|\Big|\frac{1}{2}\int_{\hat{X}}W_k{\overline W_l}\Big|\Big|$$
coincides (up to insignificant constant factor) with ${\rm det}\Im \Pi$,
where $\Pi$ is the Prym matrix (see \cite{BKN-invent}, p. 774)  corresponding to the element $(X, \Omega)$ of $Q_g(1, \dots, 1)$.
\end{Remark}

\subsection{Explicit calculation of ${\rm det}_{EKZ}\Delta^{(2)}(X, \Omega)$}
We will make use of
\begin{itemize}
\item The explicit expression for the spectral determinant of the Friedrichs extension of $\Delta^{(0)}_{|\Omega|, F}$ on the stratum $Q_g(1, \dots, 1)$
given in Theorem \ref{Leipzig}

\item Fay's version of the Belavin-Knizhnik-Manin explicit formula for the Mumford measure (more specifically, Theorem 5.8 from \cite{Fay92}).
    \end{itemize}

Assume that $g\geq 3$ (that is needed to apply Theorem 5.8 from \cite{Fay92}; the case $g=2$ can be also included but is omitted for the sake of brevity).

According to Lemmas \ref{EKZlemma} and \ref{EKZ}, the following expression
$$T=\frac{R^{(2)}_\epsilon(X, \Omega)} {\left( ({\rm det}\Im {\mathbb B})  {\rm Area}(X, |\Omega|^{(\epsilon)}\right)^{-13} R_\epsilon^{(0)}(X, \Omega)^{13}}$$
is $\epsilon$-independent. (We remind the reader that ${\rm Area}(X, |\Omega|^{(\epsilon)})={\rm Area}(X, |\Omega|)$.)

Moreover, applying Theorem 5.8 from \cite{Fay92} with $n=4$ and $\chi=I$,
one gets the relation
\begin{equation}\label{Mumford}
T=C\left|\frac{\theta\left(3K^P+\sum_{i=1}^{3g-3}{\cal A}_P(x_1+\dots+x_{3g-3})\right)\prod_{i<j}E(x_i, x_j)}
{ {\rm det}|| q_i(x_j)||\prod_{k=1}^{3g-3}{\cal C}(x_i)^{\frac{3}{g-1}}}
\right|^2\,
   \end{equation}
 where the multiplicative (g(1-g)/2)-differential ${\cal C}$ is defined in (1.17) on page 9 of \cite{Fay92};
 $x_1, \dots, x_{3g-3}$ are arbitrary points of $X$; ${\cal A}_P$ is the Abel map with base point $P$; $K^P$ is the vector of the Riemann constants and $C$ is a constant independent of $\epsilon$ and $(X, \Omega)\in Q_g(1, \dots, 1)$.

Now using explicit expression for ${\rm det}^*\Delta^{(0)}_{|\Omega|, F}$ and Proposition \ref{ravny},
one gets the following explicit expression for ${\rm det}_{EKZ}\Delta^{(2)}$:
\begin{equation}\label{detEKZ}
{\rm det}_{EKZ}\Delta^{(2)}(X, \Omega)=
\end{equation}
$$C {\rm det}||\int_X \frac{q_i \bar q_j}{|\Omega|}|| \left|\frac{\theta\left(3K^P+\sum_{i=1}^{3g-3}{\cal A}_P(x_1+\dots+x_{3g-3})\right)\prod_{i<j}E(x_i, x_j)}
{ {\rm det}\Big|\Big| q_i(x_j)\Big|\Big|\prod_{k=1}^{3g-3}{\cal C}(x_i)^{\frac{3}{g-1}}} (\tau(X, \Omega))^{13}
\right|^2\,,
$$
where $\{q_i\}$ is the Korotkin-Zograf basis of quadratic differentials on $X$ corresponding to the pair $(X, \Omega)\in Q_g(1, \dots, 1)$ and $C$ is a constant independent of $(X, \Omega)\in Q_g(1, \dots, 1)$.

\begin{Remark}  If the ${\rm det}\Delta^{(2)}_{|\Omega|}$ is understood in the sense of EKZ regularization and ${\rm det}\Delta^{(0)}_{|\Omega|}$ is understood as the $\zeta$-regularized determinant of the Friedrichs extension of the scalar conical Laplacian then the expression (\ref{Mumford0}) with conical metric $\rho$ given by $\rho=|\Omega|$, coincides with the one with any smooth conformal metric $\rho$ up to a moduli independent constant.
\end{Remark}

\end{document}